\documentclass[twoside,10pt]{article}

\usepackage{amsmath,latexsym,amssymb}

\def\demo{\noindent{\bf Proof. }}
\def\sqr#1#2{{\vcenter{\hrule height.#2pt
        \hbox{\vrule width.#2pt height#1pt \kern#1pt
                \vrule width.#2pt}
        \hrule height.#2pt}}}
\def\square{\mathchoice\sqr64\sqr64\sqr{4}3\sqr{3}3}
\def\QED{\hfill$\square$}

\newtheorem{Theorem}{Theorem}[section]
\newtheorem{Lemma}[Theorem]{Lemma}
\newtheorem{Corollary}[Theorem]{Corollary}
\newtheorem{Proposition}[Theorem]{Proposition}

\newtheorem{Remark}[Theorem]{Remark}
\newtheorem{Example}[Theorem]{Example}

\begin{document}

\title{\bf Cohen-Macaulayness of \\ special fiber rings}

\author{
\begin{tabular}{ccc}
{\normalsize\bf Alberto Corso}
 & \ &
{\normalsize\bf Laura Ghezzi} \\
\mbox{\normalsize Department of Mathematics} & \ &
\mbox{\normalsize Department of Mathematics} \\
\mbox{\normalsize University of Kentucky} & \ &
\mbox{\normalsize University of Missouri} \\
\mbox{\normalsize Lexington, KY 40506} & \ &
\mbox{\normalsize Columbia, MO 65211}  \\
\mbox{\normalsize E-mail: corso@ms.uky.edu} & \ &
\mbox{\normalsize E-mail: ghezzi@math.missouri.edu}
\end{tabular}
\\ \ \\
\begin{tabular}{ccc}
{\normalsize\bf Claudia Polini}\footnote{The last two authors
gratefully acknowledge partial support from the NSF. \newline
\indent \ \ AMS 2000 {\it Mathematics Subject Classification}.
Primary 13A30; Secondary 13B22, 13C40, 13H10.} & \ &
{\normalsize\bf Bernd Ulrich}${}^{*}$ \\
\mbox{\normalsize Department of Mathematics} & \ &
\mbox{\normalsize Department of Mathematics} \\
\mbox{\normalsize University of Notre Dame} & \ &
\mbox{\normalsize Purdue University} \\
\mbox{\normalsize Notre Dame, IN 46556} & \ &
\mbox{\normalsize West Lafayette, IN 47907}  \\
\mbox{\normalsize E-mail: cpolini@nd.edu} & \ & \mbox{\normalsize
E-mail: ulrich@math.purdue.edu}
\end{tabular}
\\ \ \\
\mbox{\normalsize\it Dedicated to Steven Kleiman on the occasion
of his sixtieth birthday} }

\date{\ }

\maketitle

\vspace{-0.5in}

\begin{abstract}
Let $(R, {\mathfrak m})$ be a Noetherian local ring and let $I$ be
an $R$-ideal. Inspired by the work of H\"ubl and Huneke, we look
for conditions that guarantee the Cohen-Macaulayness of the
special fiber ring ${\mathcal F}={\mathcal R}/{\mathfrak
m}{\mathcal R}$ of $I$, where ${\mathcal R}$ denotes the Rees
algebra of $I$. Our key idea is to require `good' intersection
properties as well as `few' homogeneous generating relations in
low degrees. In particular, if $I$ is a strongly Cohen-Macaulay
$R$-ideal with $G_{\ell}$ and the expected reduction number, we
conclude that ${\mathcal F}$ is always Cohen-Macaulay. We also
obtain a characterization of the Cohen-Macaulayness of ${\mathcal
R}/K{\mathcal R}$ for any ${\mathfrak m}$-primary ideal $K$: This
result recovers a well-known criterion of Valabrega and Valla
whenever $K=I$. Furthermore, we study the relationship among the
Cohen-Macaulay property of the special fiber ring ${\mathcal F}$
and the one of the Rees algebra ${\mathcal R}$ and the associated
graded ring ${\mathcal G}$ of $I$. Finally, we focus on the
integral closedness of ${\mathfrak m}I$. The latter question is
motivated by the theory of evolutions.
\end{abstract}


\pagestyle{myheadings}

\markright{Cohen-Macaulayness of special fiber rings} \markboth{A.
Corso, L. Ghezzi, C. Polini, and B. Ulrich}{Cohen-Macaulayness of
special fiber rings}

\section{Introduction}

\noindent Let $(R, {\mathfrak m})$ be a Noetherian local ring with
infinite residue field $k$ and let $I$ be an $R$-ideal. The {\it
Rees algebra} ${\mathcal R}$, the {\it associated graded ring}
${\mathcal G}$, and the {\it special fiber ring} ${\mathcal F}$ of
$I$,
\[
{\mathcal R}= \bigoplus_{j=0}^{\infty} I^j, \qquad {\mathcal G} =
{\mathcal R}/I{\mathcal R}, \qquad {\mathcal F} ={\mathcal
R}/{\mathfrak m}{\mathcal R}
\]
play an important role in the process of blowing up ${\rm
Spec}(R)$ along $V(I)$. For this reason these algebras are often
collectively referred to as {\it blowup algebras}.

In this article we are mostly concerned with the study of the
properties of ${\mathcal F}$. From a geometric point of view,
Proj$({\mathcal F})$ corresponds to the fiber over the closed
point of the blowup of Spec($R$) along $V(I)$. When $R$ is a
standard graded domain over $k$ and $I$ is the $R$-ideal generated
by forms $f_0, \ldots, f_n$ of the same degree, then ${\mathcal
F}$ describes the homogeneous coordinate ring of the image of the
rational map ${\rm Proj}(R) \dashrightarrow {\mathbb P}_k^n$ given
by $(f_0, \ldots, f_n)$. As a special case this construction
yields homogeneous coordinate rings of Gauss images and of secant
varieties. From a more algebraic perspective, ${\mathcal F}$
encodes $($asymptotic$)$ information about the ideal $I$: Its
Hilbert function gives the minimal number of generators of the
powers of $I$, while its Krull dimension -- referred to as the
{\it analytic spread} $\ell$ of $I$ -- coincides with the minimal
number of generators of any minimal reduction $J$ of $I$. The
analytic spread is bounded below by the height $g$ of $I$ and
bounded above by the minimal number of generators $n$ of $I$ and
by the dimension $d$ of $R$. The difference $\ell-g$ has been
dubbed {\it analytic deviation}, while $n-\ell$ is called the {\it
second analytic deviation} of $I$. The concept of a reduction $J$
of an ideal $I$ has been crucial in the study of the blowup
algebras of $I$, as $J$ carries most of the information about $I$
but, in general, with fewer generators. We say that $J$ is a {\it
reduction} of $I$ if $I^{r+1}=JI^r$ for some non-negative integer
$r$ \cite{NR}. The least such $r$ is called the reduction number
$r_J(I)$ of $I$ with respect to $J$. One then defines the {\it
reduction number} of $I$ to be the least $r_J(I)$, where $J$
varies over all minimal reductions of $I$. A reduction is said to
be {\it minimal} if it is minimal with respect to containment. We
recall that minimal reductions arise from Noether normalizations
of ${\mathcal F}$.

Special fiber rings also find an application in the theory of
evolutions. Let $k$ be a field of characteristic zero. Mazur,
inspired by the work of Wiles on semistable curves \cite{W},
introduced the notion of evolution. Let $\pi\colon S
\twoheadrightarrow T$ be an epimorphism of local $k$-algebras
essentially of finite type. We say that $S$ is an {\it evolution}
of $T$ if ${\Omega}_{S/k} \otimes_S T \cong \Omega_{T/k}$. It has
been conjectured by Mazur \cite{Mazur} that every reduced algebra
$T$ is evolutionary stable, i.e., that every evolution of $T$ is
an isomorphism. This is still an open question, although partial
results were given in \cite{EM}, \cite{HR}, \cite{Hu}. In
\cite{EM}, Eisenbud and Mazur show that for a prime ideal $I$ in a
regular local ring $(R, {\mathfrak m})$ essentially of finite type
over a field of characteristic zero, the algebra $R/I$ is
evolutionary stable if and only if $I^{(2)} \subset {\mathfrak m}
I$. As $I^{(2)}$ is contained in the integral closure
$\overline{{\mathfrak m}I}$ of ${\mathfrak m}I$, the desired
inclusion $I^{(2)} \subset {\mathfrak m}I$ holds if ${\mathfrak
m}I$ is integrally closed \cite{Hu}. H\"ubl and Huneke were the
first to use the special fiber ring ${\mathcal F}={\mathcal
R}/{\mathfrak m}{\mathcal R}$ in studying the integral closedness
of ${\mathfrak m}I$. Indeed, they observed that, if $\ell=d$,
${\mathcal R}$ is a normal domain and ${\mathfrak m}{\mathcal R}$
is unmixed, then ${\mathfrak m}{\mathcal R}$ is an integrally
closed ${\mathcal R}$-ideal, forcing all $R$-ideals ${\mathfrak
m}I^j$ to be integrally closed as well \cite{HuHu}.

\medskip

We now describe the content of the paper. Section 2, which
contains our main result, is concerned with the Cohen-Macaulay
property of ${\mathcal F}$. If $I$ is generated by part of a
system of parameters or, more generally, if $I$ has reduction
number zero $($equivalently, second analytic deviation zero$)$,
then ${\mathcal F}$ is a polynomial ring over a field, hence
Cohen-Macaulay. Thus the first interesting case is that of ideals
with reduction number one. Under this assumption, if the ring $R$
is Cohen-Macaulay, then ${\mathcal F}$ was shown to be
Cohen-Macaulay by Huneke and Sally when $I$ is ${\mathfrak
m}$-primary \cite{HS}; by Shah when $I$ is {\it equimultiple},
i.e., has analytic deviation zero \cite{shah,shah2}; by
Cortadellas and Zarzuela when $I$ has analytic deviation one and
is generically a complete intersection \cite{CZ2}.  Other cases
have been studied in \cite{MS,G,CZ3,DRV,DGH,HSw,CPV}. Recently,
H\"{u}bl and Huneke proved the Cohen-Macaulayness of ${\mathcal
F}$ for generically complete intersection ideals having analytic
deviation one, but arbitrary reduction number $r$ \cite{HuHu}.
Their key assumption is that ${\mathcal G}_{+}$ has grade $g$ and
that ${\mathcal F}$ has no homogeneous relations in degrees $\leq
r$. The work of H\"{u}bl and Huneke inspired our main result,
namely Theorem~\ref{theorem3.1}. Our main idea is that `good'
intersection properties and `few' homogeneous generating relations
in low degrees guarantee that $\ell$ general linear forms in
${\mathcal F}$ are a regular sequence, hence implying the
Cohen-Macaulayness of ${\mathcal F}$. Theorem~\ref{theorem3.1}
recovers the previous results as well as a recent one by Heinzer
and Kim \cite{HK}. We would like to thank Bill Heinzer for sharing
with us an earlier version of \cite{HK}. Our assumption on the
relations of ${\mathcal F}$ is always satisfied by ideals of
second analytic deviation one, but cannot be deleted in general,
as shown by an example due to D'Anna, Guerrieri and Heinzer $($see
Example~\ref{example3.7}$)$. Another class of ideals to which
Theorem~\ref{theorem3.1} applies are strongly Cohen-Macaulay
ideals having the `expected reduction number.' Indeed, if $I$ is a
strongly Cohen-Macaulay ideal satisfying $G_{\ell}$ and having the
`expected' reduction number $\leq \ell-g+1$, then the special
fiber ring is Cohen-Macaulay $($see
Corollary~\ref{corollary3.9}$)$.

As in \cite{HuHu}, most of our results for special fiber rings
${\mathcal F}={\mathcal R}/{\mathfrak m}{\mathcal R}$ hold more
generally for ${\mathcal R}/K{\mathcal R}$, where $K$ is an
${\mathfrak m}$-primary ideal. Under suitable assumptions we give
a characterization in terms of intersection conditions on powers
of $I$ for when ${\mathcal R}/K{\mathcal R}$ or ${\mathcal G}$ are
Cohen-Macaulay $($see Theorem~\ref{2.10} and Theorem~\ref{3.8}$)$.
For $K=I$ an ${\mathfrak m}$-primary ideal, these results
specialize to a well-known criterion by Valabrega and Valla
\cite{VV}.

In Section 3 we apply Theorem~\ref{theorem3.1} to investigate the
relationship between the Cohen-Macaulayness of ${\mathcal R}$,
${\mathcal G}$, and ${\mathcal F}$.  First of all, it is easy to
construct examples $($even for perfect ideals of height two$)$
where ${\mathcal F}$ is Cohen-Macaulay, but ${\mathcal R}$ and
${\mathcal G}$ are not. Indeed, if $I$ is a perfect ideal of
height two generated by $\ell+1$ homogeneous polynomials of the
same degree in a power series ring over a field, then ${\mathcal
F}$ is an hypersurface ring, hence Cohen-Macaulay. However,
${\mathcal R}$ $($and ${\mathcal G})$ is not Cohen-Macaulay if $I$
satisfies $G_{\ell}$, but not the `row condition' \cite{U2} $($see
Sections 2 and 3 for precise definitions$)$. Likewise, the example
by D'Anna, Guerrieri and Heinzer shows that the Cohen-Macaulayness
of ${\mathcal R}$ and ${\mathcal G}$ does not imply the
Cohen-Macaulayness of ${\mathcal F}$. Nevertheless, assuming that
$I$ satisfies $G_{\ell}$, and ${\mathcal F}$ has `few' homogeneous
generating relations in degrees $\leq {\rm max} \{ r, \ell - g
\}$, we show that if ${\mathcal G}$ is Cohen-Macaulay then so is
${\mathcal F}$ $($see Proposition~\ref{corollary3.5}$)$.
Furthermore, if $R$ is Gorenstein and $I$ is a perfect Gorenstein
ideal of height three satisfying $G_{\ell}$, then the
Cohen-Macaulayness of ${\mathcal G}$ forces $I$ to have the
expected reduction number \cite{PU}, which implies the
Cohen-Macaulayness of ${\mathcal F}$ $($see
Corollary~\ref{cor3.7}$)$. Similarly, for a perfect ideal of
height two satisfying $G_{\ell}$, ${\mathcal F}$ is Cohen-Macaulay
whenever ${\mathcal R}$ has this property $($see
Corollary~\ref{corollary3.12}$)$.

Finally, in Section 4 we identify instances when ${\mathfrak m}I$
is integrally closed. Most notably, we do this for perfect ideals
of height two and for perfect Gorenstein ideals of height three
$($see Corollary~\ref{perf2} and Corollary~\ref{gor3}$)$. To
establish both results, we use the known resolutions of the
symmetric powers of $I$ from \cite{Wey} and \cite{KU}, a criterion
given in \cite{HuHu}, and Corollary~\ref{corollary3.9}.

\medskip

\section{Cohen-Macaulayness of ${\mathcal R}/K{\mathcal R}$}

\noindent Our main result is Theorem~\ref{theorem3.1} below. It
fully generalizes \cite[2.1]{HuHu}. To prove it we need a
preparatory lemma.

\begin{Lemma}\label{lemma2.1}
Let $R$ be a Noetherian ring, let $I$ be an $R$-ideal, let $a_1,
\ldots, a_{\ell}$ be elements generating a reduction $J$ of $I$,
and let $t$ be a fixed positive integer with $t \geq r_J(I)$. If
\[
[(a_{\nu_1}, \ldots, a_{\nu_i})I^{t-1} \colon a_{\nu_{i+1}}] \cap
I^j = (a_{\nu_1}, \ldots, a_{\nu_i})I^{j-1}
\]
holds for every $\{ \nu_1, \ldots, \nu_{i+1} \} \subset \{ 1,
\ldots, \ell\}$ and for $j=t$, then it holds for every $j \geq t$.
\end{Lemma}
\demo We proceed by induction on $j \geq t$. The case $j=t$ is
satisfied by assumption, thus let $j \geq t+1$. We now use
decreasing induction on $i \leq \ell-1$. Write $L=[(a_{\nu_1},
\ldots, a_{\nu_i})I^{t-1} \colon a_{\nu_{i+1}}] \cap I^j$ and
choose $\nu_{i+2}=\nu_{i+1}$ if $i=\ell-1$ and $\nu_{i+2} \in \{
1, \ldots, \ell \} \setminus \{ \nu_1, \ldots, \nu_{i+1} \}$ if
$i<\ell-1$. One has
\begin{eqnarray*}
L & \subset & [(a_{\nu_1}, \ldots, a_{\nu_i},
a_{\nu_{i+2}})I^{t-1} \colon a_{\nu_{i+1}}] \cap I^j \\
& = & (a_{\nu_1}, \ldots, a_{\nu_i}, a_{\nu_{i+2}})I^{j-1}.
\end{eqnarray*}
Indeed, for $i=\ell-1$ the last equality holds since $j \geq
r_J(I)+1$ while for $i<\ell-1$ it holds by decreasing induction on
$i$. Since $(a_{\nu_1}, \ldots, a_{\nu_i})I^{j-1} \subset L$ we
obtain
\begin{eqnarray*}
L & = & (a_{\nu_1}, \ldots, a_{\nu_i},
a_{\nu_{i+2}})I^{j-1} \cap L \\
& = & (a_{\nu_1}, \ldots, a_{\nu_i})I^{j-1}
+ (a_{\nu_{i+2}} I^{j-1} \cap L) \\
& = & (a_{\nu_1}, \ldots, a_{\nu_i})I^{j-1}+a_{\nu_{i+2}}[(L
\colon a_{\nu_{i+2}}) \cap I^{j-1}].
\end{eqnarray*}
By assumption $L \subset (a_{\nu_1}, \ldots, a_{\nu_i})I^{t-1}$.
Hence
\[
(L \colon a_{\nu_{i+2}}) \cap I^{j-1} \subset [(a_{\nu_1}, \ldots,
a_{\nu_i})I^{t-1} \colon a_{\nu_{i+2}}] \cap I^{j-1} = (a_{\nu_1},
\ldots, a_{\nu_i})I^{j-2},
\]
where the last equality holds by induction on $j$. Thus
\begin{eqnarray*}
L & \subset & (a_{\nu_1}, \ldots, a_{\nu_i})I^{j-1}
+ a_{\nu_{i+2}}(a_{\nu_1}, \ldots, a_{\nu_i})I^{j-2} \\
& = & (a_{\nu_1}, \ldots, a_{\nu_i})I^{j-1},
\end{eqnarray*}
as desired. \QED

\bigskip

Let $R$ be a Noetherian local ring and let $I$ be an $R$-ideal
minimally generated by $f_1, \ldots, f_n$. Writing ${\mathcal A}$
for the kernel of the homogeneous epimorphism $R[T_1, \ldots, T_n]
\twoheadrightarrow {\mathcal R}$ of $R$-algebras which sends $T_i$
to $f_i$, one has ${\mathcal R} \cong R[T_1, \ldots,
T_n]/{\mathcal A}$. Thus for any $R$-ideal $K$, ${\mathcal
R}/K{\mathcal R} \cong R[T_1,\dots,T_n]/(K,{\mathcal A} )\cong
(R/K)[T_1,\dots,T_n]/Q$, where $Q$ denotes the image in
$(R/K)[T_1, \ldots, T_n]$ of ${\mathcal A}$. By `homogeneous
generating relations' of ${\mathcal R}/K{\mathcal R}$ we will mean
elements forming part of a homogeneous minimal generating set of
the ideal $Q$.

In the sequel we will denote the minimal number of generators
function by $\mu$.

\begin{Theorem}\label{theorem3.1}
Let $(R, {\mathfrak m})$ be a Noetherian local ring with infinite
residue field, let $I$ be an $R$-ideal with analytic spread
$\ell$, minimal number of generators $n$, reduction number $r$,
and let $t$ be a positive integer with $t\geq r$. Let $a_1,\dots,
a_{\ell}$ be general elements in $I$ and let ${\mathfrak
a}_i=(a_1,\dots,a_i)$ for $0 \leq i \leq \ell$. Let $K$ be an
${\mathfrak m}$-primary ideal. Assume that
\begin{itemize}
\item[$({\it a})$]
$({\mathfrak a}_iI^{t-1} \colon a_{i+1})\cap I^t = {\mathfrak
a}_iI^{t-1}$ whenever $0\leq i\leq \ell-1$;

\item[$({\it b})$]
\begin{itemize}
\item[$({\it i})$]
if $K={\mathfrak m}$ and $n \geq \ell+2$ then ${\mathcal
R}/K{\mathcal R}={\mathcal F}$ has at most two homogeneous
generating relations in degrees $\leq t$;

\item[$({\it ii})$]
if $K={\mathfrak m}$ and $n=\ell+1$ then ${\mathcal R}/K{\mathcal
R }={\mathcal F}$ has at most one homogeneous generating relation
in degrees $\leq t$;

\item[$({\it iii})$]
if $K \not= {\mathfrak m}$ then ${\mathcal R}/K{\mathcal R}$ has
no homogeneous relations in degrees $\leq t$.
\end{itemize}
\end{itemize}
Then ${\mathcal R}/K{\mathcal R}$ is Cohen-Macaulay.
\end{Theorem}
\demo Set $S=(R/K)[T_1, \ldots, T_n]$, ${\mathcal R}/K{\mathcal
R}=S/Q$, let $'$ denote images under the homomorphism ${\mathcal
R} \longrightarrow {\mathcal R}/K{\mathcal R}$, and write $Q_{\leq
t}$ for the $S$-ideal generated by the forms in $Q$ of degrees at
most $t$. If $K \not= {\mathfrak m}$ then $S/Q_{\leq t}=S$ by
assumption $({\it b})({\it iii})$. If $K={\mathfrak m}$ and $n
\geq \ell + 2$, then $\mu (Q_{\leq t}) \leq 2$ according to
assumption $({\it b})({\it i})$. Hence ${\rm proj\,dim}_S \,
Q_{\leq t} \leq 1$ and so ${\rm depth}\, S/Q_{\leq t}\geq n-2\geq
\ell$. Similarly if $K={\mathfrak m}$ and $n = \ell + 1$, then
$\mu (Q_{\leq t}) \leq 1$ by assumption $({\it b})({\it ii})$, and
${\rm depth}\, S/Q_{\leq t}\geq n-1=\ell$. In any case we may
assume that the images $a_1',\dots,a_{\ell}'$ of $a_1,\dots,
a_{\ell}$ in $I/KI = [S/Q]_1$ form a regular sequence on the ring
$S/Q_{\leq t}$.

Next we claim that
\begin{equation}\label{eq1}
{\mathfrak a}_iI^j\cap KI^{j+1}= {\mathfrak a}_iKI^j \ \ {\rm
whenever} \ \ 0 \leq i \leq \ell \ \ {\rm and} \ \ j\geq 0.
\end{equation}
The proof is by induction on $j$.

\noindent First assume $j\leq t-1$. Let
$\lambda_1a_1+\dots+\lambda_ia_i\in KI^{j+1}=[K{\mathcal
R}]_{j+1}$ with $\lambda_1,\dots,\lambda_i$ elements of
$I^j=[{\mathcal R}]_j$. One has
$\lambda_1'a_1'+\dots+\lambda_i'a_i'=0$ in $S/Q_{\leq t}$ since
$j+1\leq t$. As $a_1',\dots, a_i'$ form a regular sequence in this
ring, there is an alternating $i$ by $i$ matrix $A$ with entries
in $[{\mathcal R}]_{j-1}$ so that
$[\lambda_1',\dots,\lambda_i']=[a_1',\dots,a_i']A'$. Hence
$[\lambda_1,\dots,\lambda_i] \equiv [a_1,\dots,a_i]A$ mod
$KI^{j}$. As $[a_1,\dots,a_i]A[a_1,\dots,a_i]^{\rm tr}=0$, it
follows that $\lambda_1a_1+\dots+\lambda_ia_i\in {\mathfrak
a}_iKI^j$.

\noindent Next suppose $j\geq t$. We use decreasing induction on
$i$, the case $i=\ell$ being clear since $j\geq r$. If $i<\ell$
then \newline
\[
\begin{array}{rcll} \vspace{.3cm}
{\mathfrak a}_iI^j\cap KI^{j+1} & = & {\mathfrak a}_iI^j \cap
{\mathfrak a}_{i+1}I^j \cap KI^{j+1} & \\ \vspace{.3cm} & = &
{\mathfrak a}_iI^j \cap {\mathfrak a}_{i+1}KI^j & \text{by
induction on} \ i
\\ \vspace{.3cm}
& = & {\mathfrak a}_iI^j \cap ({\mathfrak a}_i KI^j + a_{i+1}KI^j)
& \\ \vspace{.3cm} & = & {\mathfrak a}_iKI^j +({\mathfrak a}_iI^j
\cap a_{i+1}KI^j) &
\\ \vspace{.3cm}
& = & {\mathfrak a}_i KI^{j}+a_{i+1}[({\mathfrak
a}_iI^j:a_{i+1})\cap KI^j] & \\ \vspace{.1cm} & = & {\mathfrak
a}_iKI^{j}+ a_{i+1}({\mathfrak a}_iI^{j-1}\cap KI^j) & \text{by
assumption } ({\it a}) \\ \vspace{.3cm} & & & \text{and
Lemma~\ref{lemma2.1}}
\\ \vspace{.3cm}
& = & {\mathfrak a}_iKI^{j}+a_{i+1}{\mathfrak a}_iKI^{j-1} &
\text{by induction on } j \\ \vspace{.3cm} & = & {\mathfrak
a}_iKI^{j}. &
\end{array}
\]
This completes the proof of $(\ref{eq1})$.

Finally we prove that ${\mathcal R}/K{\mathcal R}$ is
Cohen-Macaulay. Recall that $a_1',\dots, a_{\ell}'$ form a regular
sequence on $S/Q_{\leq t}$. Hence it suffices to show that for
$0\leq i\leq \ell-1$ and $j\geq t$,
\[
[(a_1',\dots, a_i') :_{{\mathcal R}/K{\mathcal R}}
a_{i+1}']_j=[(a_1',\dots, a_i')]_j,
\]
or equivalently,
\[
[({\mathfrak a}_iI^j+KI^{j+1}) :_R a_{i+1}] \cap I^j \subset
{\mathfrak a}_iI^{j-1}+KI^j.
\]
Notice that
\begin{eqnarray*}
({\mathfrak a}_iI^j+KI^{j+1}) \cap {\mathfrak a}_{i+1}I^j & = &
{\mathfrak a}_iI^j + ({\mathfrak a}_{i+1}I^j \cap KI^{j+1}) \\
& = & {\mathfrak a}_iI^j+{\mathfrak a}_{i+1}KI^j \qquad \qquad
\quad {\rm by}
\ (\ref{eq1}) \\
& = & {\mathfrak a}_iI^j+a_{i+1}KI^j.
\end{eqnarray*}
Thus
\begin{eqnarray*}
[({\mathfrak a}_iI^j + KI^{j+1}) \colon a_{i+1}]\cap I^j & \!\!\!
\subset \!\!\! & [(({\mathfrak a}_iI^j + KI^{j+1}) \cap {\mathfrak
a}_{i+1}I^j)
\colon a_{i+1}]\cap I^j \\
& \!\!\! = \!\!\! & [({\mathfrak a}_iI^j+a_{i+1}KI^j) \colon
a_{i+1}] \cap I^j \hspace{0.4cm}
{\rm by \ the \ above} \\
& \!\!\! = \!\!\! & [({\mathfrak a}_iI^j \colon a_{i+1}) + KI^j] \cap I^j\\
& \!\!\! = \!\!\! & ({\mathfrak a}_iI^j \colon a_{i+1})\cap I^j+KI^j \\
& \!\!\! = \!\!\! & {\mathfrak a}_iI^{j-1}+KI^j,
\end{eqnarray*}
where the last equality follows from assumption $({\it a})$ and
Lemma~\ref{lemma2.1}. \QED

\bigskip

One often does not have to require assumption $({\it b})$ in
Theorem~\ref{theorem3.1} if the ideal $I$ has second analytic
deviation one. Notice that in this case ${\mathcal F}$ is
Cohen-Macaulay if and only if ${\mathcal F}$ is a hypersurface
ring.

\begin{Corollary}\label{corollary3.2}
Let $R$ be a Noetherian local ring with infinite residue field,
and let $I$ be an $R$-ideal with analytic spread $\ell$, minimal
number of generators $\ell+1$ and reduction number $r$. Let
$a_1,\dots, a_{\ell}$ be general elements in $I$ and let
${\mathfrak a}_i=(a_1,\dots,a_i)$ for $0 \leq i \leq \ell$. Assume
that $({\mathfrak a}_iI^{r-1} \colon a_{i+1}) \cap I^r =
{\mathfrak a}_iI^{r-1}$ whenever $0 \leq i \leq \ell-1$. Then
${\mathcal F}$ is Cohen-Macaulay.
\end{Corollary}
\demo As $\mu(I) = \ell+1$, ${\mathcal F}$ has no homogeneous
relations in degrees $\leq r$. Thus our assertion follows from
Theorem~\ref{theorem3.1} with $t=r$. \QED

\begin{Remark}
{\rm Theorem~\ref{theorem3.1} recovers \cite[3.2]{CZ3} and
\cite[2.1]{HuHu}. Indeed, the assumptions of \cite[3.2]{CZ3} and
\cite[2.1]{HuHu} imply $({\mathfrak a}_i \colon a_{i+1}) \cap I^j
= {\mathfrak a}_iI^{j-1}$ whenever $0 \leq i \leq \ell-1$ and
$j\geq 1$. }
\end{Remark}
\demo In the case of H\"ubl and Huneke, we have that $\ell=g+1$,
where $g$ is the height of $I$. Moreover, ${\rm ht}\,(({\mathfrak
a}_g \colon I ), a_{g+1}) \geq g+1$, $a_1, \ldots, a_g$ form a
regular sequence on $R$, and their images in $[{\mathcal
G}]_1=I/I^2$ form a regular sequence on ${\mathcal G}$. As $i \leq
\ell-1 = g$, the first and second condition give $({\mathfrak a}_i
\colon a_{i+1}) \cap I = {\mathfrak a}_i$, whereas the third one
implies ${\mathfrak a}_i \cap I^j ={\mathfrak a}_iI^{j-1}$ for
every $j \geq 1$ by \cite[2.7]{VV}.

In the case of Cortadellas and Zarzuela the intersection
properties follow from \cite[2.5({\it ii})]{CZ3}. \QED

\medskip

\begin{Remark}
{\rm Corollary~\ref{corollary3.2} also recovers the following
result of Heinzer and Kim \cite[5.6]{HK}: Let $R$ be a Noetherian
local ring with infinite residue field and $I$ an $R$-ideal of
grade $g$ with analytic spread $g$ and minimal number of
generators $g+1$. If ${\rm grade}\, {\mathcal G}_{+} \geq g-1$
then ${\mathcal F}$ is Cohen-Macaulay. In particular ${\mathcal
F}$ is a hypersurface ring. }
\end{Remark}
\demo We use the notation of Corollary~\ref{corollary3.2} with
$\ell=g$. Notice that ${\mathfrak a}_i \colon a_{i+1} = {\mathfrak
a}_i$ for $0 \leq i \leq g-1$, since $a_1, \ldots, a_g$ form a
regular sequence; moreover, ${\mathfrak a}_i \cap I^j ={\mathfrak
a}_iI^{j-1}$ for $0 \leq i \leq g-1$ and $j \geq 1$ since ${\rm
grade}\, {\mathcal G}_{+} \geq g-1$ \cite[2.7]{VV}. Now the
assertion follows from Corollary~\ref{corollary3.2}. \QED

\bigskip

The work of H\"ubl and Swanson \cite{HSw} provides classes of
defining ideals of monomial space curves for which
Theorem~\ref{theorem3.1} implies the Cohen-Macaulayness of
${\mathcal F}$ in the presence of two homogeneous generating
relations of degrees $\leq r$.

Next, we use an example by D'Anna, Guerrieri and Heinzer
\cite[2.3]{DGH} to point out that in general the assumption
`${\mathcal F}$ has at most two homogeneous generating relations
in degrees $\leq r$' in Theorem~\ref{theorem3.1} cannot be removed
or weakened.

\begin{Example}\label{example3.7}
{\rm Let $R=k[\![t^6,t^{11},t^{15},t^{31}]\!]$, where $k$ is an
infinite field, and let $I=(t^6,t^{11},t^{31})$. The ideal $I$ has
height 1, analytic spread 1, minimal number of generators $3$ and
reduction number 2. Assumption $({\it a})$ of
Theorem~\ref{theorem3.1} is satisfied, but condition $({\it
b})({\it ii})$ does not hold. Indeed, one has
\[
{\mathcal F}= k[T_1,T_2,T_3]/(T_2^3,T_1T_3,T_2T_3,T_3^2);
\]
so ${\mathcal F}$ has 3 generating relations in degree 2. Notice
that ${\mathcal F}$ is not Cohen-Macaulay. }
\end{Example}

In contrast to Theorem~\ref{theorem3.1}, the following criterion
for the Cohen-Macau\-layness of ${\mathcal R}/K{\mathcal R}$ does
not refer to the degrees of the defining equations. Instead, we
impose conditions on certain intersections that specialize to the
assumptions of Valabrega and Valla \cite[2.7]{VV} whenever $K=I$.

\begin{Theorem}\label{2.10}
Let $(R, {\mathfrak m})$ be a Noetherian local ring with infinite
residue field, and let $I$ be an $R$-ideal with analytic spread
$\ell$ and reduction number $r$. Let $a_1,\dots, a_{\ell}$ be
general elements in $I$ and set $J=(a_1, \ldots, a_{\ell})$. Let
${\mathfrak a}_i=(a_1,\dots,a_i)$ for $0 \leq i \leq \ell$, and
assume that
\[
({\mathfrak a}_iI^{j-1} \colon a_{i+1})\cap I^j = {\mathfrak
a}_iI^{j-1}
\]
whenever $0\leq i \leq \ell-1$ and $1 \leq j \leq \max\{1, r\}$.
Let $K$ be an ${\mathfrak m}$-primary ideal. Then the following
are equivalent:
\begin{itemize}
\item[$({\it i})$]
${\mathcal R}/K{\mathcal R}$ is Cohen-Macaulay;

\item[$({\it ii})$]
$J/KJ$ is a free $R/K$-module and $JI^{j-1} \cap KI^j=JKI^{j-1}$
whenever $1 \leq j \leq r$.
\end{itemize}
\end{Theorem}
\demo We first show that $({\it i})$ implies $({\it ii})$. Notice
that the images $a_1', \ldots, a_{\ell}'$ of $a_1, \ldots,
a_{\ell}$ in $[{\mathcal R}/K{\mathcal R}]_1 = I/KI$ form a
regular sequence on ${\mathcal R}/K{\mathcal R}$. Thus $a_1',
\ldots, a_{\ell}'$ are linearly independent over $[{\mathcal
R}/K{\mathcal R}]_0=R/K$. In particular, the images of $a_1,
\ldots, a_{\ell}$ in $J/KJ$ are linearly independent  over $R/K$
and hence form a basis of $J/KJ$. Again since $a_1', \ldots,
a_{\ell}'$ form a regular sequence on ${\mathcal R}/K{\mathcal
R}$, the equalities $JI^{j-1} \cap KI^j=JKI^{j-1}$ for $j \geq 1$
follow as in the proof of Theorem~\ref{theorem3.1} $($see the
proof of claim $($\ref{eq1}$)$ in the case $j \leq t-1)$.

Next we show that $({\it ii})$ implies $({\it i})$. Our assumption
and Lemma~\ref{lemma2.1} yield
\begin{equation}\label{eq2}
({\mathfrak a}_iI^{j-1} \colon a_{i+1}) \cap I^j ={\mathfrak a}_i
I^{j-1} \qquad {\rm whenever} \ 0 \leq i \leq \ell-1 \ {\rm and} \
j \geq 1.
\end{equation}
We claim that
\begin{equation}\label{eq3}
{\mathfrak a}_iI^{j-1} \cap KI^j = {\mathfrak a}_i K I^{j-1} \quad
\! {\rm whenever} \ 0 \leq i \leq \ell \ {\rm and} \ j \geq 1.
\end{equation}
By the freeness of $J/KJ$ one has ${\mathfrak a}_i \cap KJ =
{\mathfrak a}_i K$, whereas the second assumption in $({\it ii})$
gives $J \cap KI = KJ$. Thus ${\mathfrak a}_i \cap KI = {\mathfrak
a}_i \cap J \cap KI = {\mathfrak a}_i \cap KJ = {\mathfrak a}_i
K$, which yields $($\ref{eq3}$)$ for $j=1$. On the other hand,
$($\ref{eq3}$)$ holds for $i=\ell$ by the second assumption in
$({\it ii})$. Now $($\ref{eq3}$)$ follows as in the proof of
Theorem~\ref{theorem3.1} $($see the proof of claim $(1)$ in the
case $j \geq t)$, by increasing induction on $j$ and decreasing
induction on $i$ using $($\ref{eq2}$)$. Finally, again the proof
of Theorem~\ref{theorem3.1} shows that $($\ref{eq2}$)$ and
$($\ref{eq3}$)$ imply the Cohen-Macaulayness of ${\mathcal
R}/K{\mathcal R}$. \QED

\bigskip

We now study several instances where assumption $({\it a})$ in
Theorem~\ref{theorem3.1} as well as the assumption in
Theorem~\ref{2.10} are automatically satisfied.

\medskip

Let $R$ be a Noetherian local ring, $I$ an $R$-ideal of height
$g$, and $s$ an integer. Recall that $I$ satisfies condition $G_s$
if $\mu(I_{{\mathfrak p}}) \leq {\rm dim}\, R_{{\mathfrak p}}$ for
every ${\mathfrak p} \in V(I)$ with ${\rm dim}\, R_{\mathfrak p}
\leq s-1$. A proper $R$-ideal $H$ is called an {\it $s$-residual
intersection} of $I$, if there exists an $s$-generated ideal
${\mathfrak a} \subset I$ so that $H = {\mathfrak a} \colon I$ and
${\rm ht}\, H \geq s \geq g$. If in addition ${\rm ht}\, (I + H)
\geq s+1$, then $H$ is said to be a {\it geometric $s$-residual
intersection} of $I$. We say that $I$ satisfies $AN_s$
$(AN_s^{-}$, respectively$)$ if $R/H$ is Cohen-Macaulay for every
$i$-residual intersection $($geometric $i$-residual intersection,
respectively$)$ $H$ of $I$ and every $i \leq s$.

\begin{Lemma}\label{lemma2.7}
Let $R$ be a local Cohen-Macaulay ring with infinite residue field
and let $I$ be an $R$-ideal. Let $s, t$ be integers. Let $J$ be a
reduction of $I$ generated by $s$ elements with ${\rm ht}\, J
\colon I \geq s$. Assume that $I$ satisfies $G_s$,
$AN_{s-t-1}^{-}$ and that for every ${\mathfrak p} \in V(I)$,
${\rm depth}\, (R/I^j)_{{\mathfrak p}} \geq \min \{ {\rm dim}\,
R_{{\mathfrak p}}-s+t-j, t-j \}$ whenever $1 \leq j \leq t-1$. Let
$a_1, \ldots, a_s$ be general elements in $J$ and set ${\mathfrak
a}_i =(a_1, \ldots, a_i)$ for $0 \leq i \leq s$. Then $({\mathfrak
a}_i \colon a_{i+1}) \cap I^j={\mathfrak a}_iI^{j-1}$ whenever $0
\leq i \leq s-1$ and $\max\{ 1, i-s+t+1\} \leq j \leq t$.
\end{Lemma}
\demo See the proof of \cite[2.4]{L}. \QED

\medskip

\begin{Corollary}\label{corollary3.6}
Let $(R, {\mathfrak m})$ be a local Cohen-Macaulay ring with
infinite resi\-due field, let $I$ be an $R$-ideal with analytic
spread $\ell$, minimal number of generators $n$ and reduction
number $r$, and let $t$ be a positive integer with $t \geq r$.
Assume that $I$ satisfies $G_{\ell}$, $AN^-_{\ell-t-1}$ and that
for every ${\mathfrak p} \in V(I)$, ${\rm depth}\,
(R/I^j)_{{\mathfrak p}} \geq \min\{\dim R_{{\mathfrak
p}}-\ell+t-j, t-j\}$ whenever $1\leq j\leq t-1$. Let $K$ be an
${\mathfrak m}$-primary ideal. Assume that
\begin{itemize}
\item[$({\it i})$]
if $K={\mathfrak m}$ and $n \geq \ell+2$ then ${\mathcal
R}/K{\mathcal R}={\mathcal F}$ has at most two homogeneous
generating relation in degrees $\leq t$;

\item[$({\it ii})$]
if $K={\mathfrak m}$ and $n=\ell+1$ then ${\mathcal R}/K{\mathcal
R }={\mathcal F}$ has at most one homogeneous generating relation
in degrees $\leq t$;

\item[$({\it iii})$]
if $K \not= {\mathfrak m}$ then ${\mathcal R}/K{\mathcal R}$ has
no homogeneous relations in degrees $\leq t$.
\end{itemize}
Then ${\mathcal R}/K{\mathcal R}$ is Cohen-Macaulay.
\end{Corollary}
\demo We use the notation of Theorem~\ref{theorem3.1}. By
Lemma~\ref{lemma2.7} we conclude that $({\mathfrak a}_i \colon
a_{i+1})\cap I^j = {\mathfrak a}_iI^{j-1}$ and, a fortiori,
$({\mathfrak a}_iI^{j-1} \colon a_{i+1})\cap I^j = {\mathfrak
a}_iI^{j-1}$ whenever $0\leq i\leq \ell-1$ and
$\max\{1,i-\ell+t+1\} \leq j \leq t$. Our assertion then follows
from Theorem~\ref{theorem3.1}. \QED

\medskip

\begin{Corollary}\label{corollary3.7}
Let $(R, {\mathfrak m})$ be a local Cohen-Macaulay ring with
infinite resi\-due field, and let $I$ be an $R$-ideal of height
$g$ with analytic spread $\ell$, minimal number of generators $n$
and reduction number $r$. Assume that $I$ satisfies $G_{\ell}$.
Further suppose that either $t=\ell-g+1 \geq r$ and ${\rm depth}\,
R/I^j \geq {\rm dim}\, R/I -j+1$ whenever $1 \leq j \leq \ell-g$,
or else $t=\ell-g \geq r$ and ${\rm depth}\, R/I^j \geq {\rm
dim}\, R/I -j$ whenever $1 \leq j \leq \ell-g-1$. Let $K$ be an
${\mathfrak m}$-primary ideal containing $I$. Assume that
\begin{itemize}
\item[$({\it i})$]
if $K={\mathfrak m}$ and $n \geq \ell+2$ then ${\mathcal
R}/K{\mathcal R}={\mathcal F}$ has at most two homogeneous
generating relations in degrees $\leq t$;

\item[$({\it ii})$]
if $K={\mathfrak m}$ and $n=\ell+1$ then ${\mathcal R}/K{\mathcal
R }={\mathcal F}$ has at most one homogeneous generating relation
in degrees $\leq t$;

\item[$({\it iii})$]
if $K \not= {\mathfrak m}$ then ${\mathcal R}/K{\mathcal R}$ has
no homogeneous relations in degrees $\leq t$.
\end{itemize}
Then ${\mathcal R}/K{\mathcal R}$ is Cohen-Macaulay.
\end{Corollary}
\demo We may assume $t>0$ since otherwise $I$ is a complete
intersection. The assertion now follows from
Corollary~\ref{corollary3.6}. Indeed the assumption
$AN_{\ell-t-1}^{-}$ of that corollary is vacuous since $\ell-t-1
\leq g-1$, whereas the local depth condition follows from the
global one by Ischebeck's Lemma, for instance. \QED

\bigskip

In the presence of the condition $G_{\ell}$, the depth assumptions
of Corollary~\ref{corollary3.7} are satisfied by {\it strongly
Cohen-Macaulay ideals}, i.e., ideals whose Koszul homology modules
are Cohen-Macaulay. This condition always holds if $I$ is a
Cohen-Macaulay almost complete intersection or a Cohen-Macaulay
deviation two ideal of a Gorenstein ring \cite[p.{\,}259]{AH}. It
is also satisfied for any ideal in the linkage class of a complete
intersection \cite[1.11]{H1}: Standard examples include perfect
ideals of height two and perfect Gorenstein ideals of height
three.

\medskip

Continuing the discussion preceding Theorem~\ref{theorem3.1}, let
$R$ be a Noetherian local ring, let $I$ be an $R$-ideal minimally
generated by $f_1, \ldots, f_n$, and let $K$ be any $R$-ideal.
There are homogeneous epimorphisms of graded $R$-algebras,
\[
B=R[T_1, \ldots, T_n] \twoheadrightarrow {\rm Sym}(I) =
\displaystyle\bigoplus_{j=0}^{\infty} S_j(I) \twoheadrightarrow
{\mathcal R},
\]
where the first map sends $T_i$ to $f_i$ and the
second map is the natural one. They yield presentations
\[
{\mathcal R} \cong B/{\mathcal A}, \quad {\rm Sym}(I) \cong
B/{\mathcal A}_{\leq 1}, \quad {\mathcal R}/K{\mathcal R} \cong
(R/K)[T_1, \ldots, T_n]/Q.
\]
Here ${\mathcal A}_{\leq 1}$ denotes the $B$-ideal generated by
$[{\mathcal A}]_1$, and $Q$ is the image of ${\mathcal A}$ in
$(R/K)[T_1, \ldots, T_n]$. One has ${\mathcal A}_{\leq 1} \subset
{\rm Fitt}_{n-1}(I)B$. It follows that if ${\rm Fitt}_{n-1}(I)
\subset K$ and $S_j(I) \cong I^j$ for some $j$, then $[Q]_j=0$.

\begin{Corollary}\label{corollary3.9}
Let $(R, {\mathfrak m})$ be a local Cohen-Macaulay ring with
infinite residue field, and let $I$ be a strongly Cohen-Macaulay
$R$-ideal of height $g$ with analytic spread $\ell$, minimal
number of generators $n$ and reduction number $\leq \ell-g+1$.
Assume that $I$ satisfies $G_{\ell}$. Then ${\mathcal
R}/K{\mathcal R}$ is Cohen-Macaulay for every ${\mathfrak
m}$-primary ideal $K$ containing ${\rm Fitt}_{n-1}(I)$. In
particular, ${\mathcal F}$ is Cohen-Macaulay.
\end{Corollary}
\demo We apply Corollary~\ref{corollary3.7} with $t=\ell-g+1$.
Indeed, \cite[the proofs of  5.1 and 4.6]{HSV} shows that ${\rm
depth}\, R/I^j \geq {\rm dim}\, R/I-j+1$ and $S_j(I) \cong I^j$
for $1 \leq j \leq \ell-g+1$. Since ${\rm Fitt}_{n-1}(I) \subset
K$ the last isomorphisms imply that ${\mathcal R}/K{\mathcal R}$
has no homogeneous relations in degrees $\leq \ell-g+1$. \QED

\bigskip

The next example shows that Corollary~\ref{corollary3.9} is no
longer true even for perfect ideals of height two with second
analytic deviation one, if the reduction number is not the
`expected' one.

\begin{Example}\label{example 3.10}
{\rm Let $R=k[\![x,y,z]\!]$, with $k$ an infinite field, and let
$I$ be the $R$-ideal generated by the $3$ by $3$ minors of the
matrix
\[
\varphi=
\begin{pmatrix}x^3 & 0 & 0\cr
       y^2 & 0 & yz\cr
       0 &y^2 & z^2\cr
       0 & z^2 & x^2\cr
       \end{pmatrix}.
\]
One has that $I$ is a perfect ideal of height $2$, analytic spread
$3$, and reduction number $5$. Also, $I$ satisfies $G_3$. However
${\mathcal R}/K{\mathcal R}$ is not Cohen-Macaulay for any
$R$-ideal $K$ containing ${\rm Fitt}_3(I) =(x^2, y^2, z^2, yz)$.
Indeed, writing $A=R/K$ we have that
\[
{\mathcal R}/(K, T_1, T_3){\mathcal R} \cong A[T_2,
T_4]/(zT_2^2T_4, yT_2^2T_4, xT_2^3T_4, T_2^5T_4, T_2^4T_4^2)
\]
has depth zero and dimension one. As $T_1$ and $T_3$ form part of
a system of parameters of ${\mathcal R}/K{\mathcal R}$ it follows
that ${\mathcal R}/K{\mathcal R}$ cannot be Cohen-Macaulay. In
particular, the special fiber ring ${\mathcal F}$ is not
Cohen-Macaulay.}
\end{Example}

\medskip

\section{The Cohen-Macaulayness of ${\mathcal R}/K{\mathcal R}$
versus the Cohen-Macaulayness of ${\mathcal R}$ and ${\mathcal G}$
\hspace{5.15in} }

Our goal in this section is to relate the Cohen-Macaulayness of
${\mathcal R}$ and ${\mathcal G}$ to the Cohen-Macaulayness of
${\mathcal R}/K{\mathcal R}$ and vice versa.

\begin{Proposition}\label{corollary3.5}
Let $(R, {\mathfrak m})$ be a Noetherian local ring with infinite
residue field, and let $I$ be an $R$-ideal of height $g$ with
analytic spread $\ell$, minimal number of generators $n$ and
reduction number $r$. Assume that $I$ satisfies $G_{\ell}$. Let
$K$ be an ${\mathfrak m}$-primary $R$-ideal containing $I$.
Suppose that
\begin{itemize}
\item[$({\it i})$]
if $K={\mathfrak m}$ and $n \geq \ell+2$ then ${\mathcal
R}/K{\mathcal R}={\mathcal F}$ has at most two homogeneous
generating relations in degrees $\leq \max \{r, \ell-g\}$;

\item[$({\it ii})$]
if $K={\mathfrak m}$ and $n=\ell+1$ then ${\mathcal R}/K{\mathcal
R }={\mathcal F}$ has at most one homogeneous generating relation
in degrees $\leq \ell-g$;

\item[$({\it iii})$]
if $K \not= {\mathfrak m}$ then ${\mathcal R}/K{\mathcal R}$ has
no homogeneous relations in degrees $\leq \max\{r,\ell-g\}$.
\end{itemize}
If $\mathcal G$ is Cohen-Macaulay, then ${\mathcal R}/K{\mathcal
R}$ is Cohen-Macaulay.
\end{Proposition}
\demo We use the notation of Theorem~\ref{theorem3.1}. From
\cite[2.2]{J} we conclude that $({\mathfrak a}_i:a_{i+1}) \cap
I^j= {\mathfrak a}_iI^{j-1}$ whenever $0\leq i\leq \ell-1$ and
$j\geq i-g+1$. Notice that in the setting of $({\it ii})$,
${\mathcal F}$ has no homogeneous relations in degrees $\leq r$.
Furthermore we may assume $\max\{ r, \ell-g \} > 0$ since
otherwise $I$ is a complete intersection. Now our assertion
follows from Theorem~\ref{theorem3.1} with $t=\max\{ r, \ell-g
\}$. \QED

\medskip

\begin{Corollary}\label{corollaryfiber}
Let $R$ be a Noetherian local ring with infinite residue field,
and let $I$ be an $R$-ideal of height $g$ with analytic spread
$\ell$ and reduction number $r\geq \ell-g$. Assume that $I$
satisfies $G_{\ell}$ and that ${\mathcal F}$ has at most two
homogeneous generating relations in degrees $\leq r$. If
${\mathcal G}$ is Cohen-Macaulay, then $\mathcal F$ is
Cohen-Macaulay.
\end{Corollary}
\demo Again recall that if $n=\ell+1$ then ${\mathcal F}$ has no
homogeneous generating relations in degrees $\leq r$. \QED

\medskip

Example~\ref{example3.7} shows once more that the assumption
`${\mathcal F}$ has at most two generating relations in degrees
$\leq r$' cannot be weakened in Corollary~\ref{corollaryfiber}.
That is, ${\mathcal F}$ may fail to be Cohen-Macaulay even if
${\mathcal G}$ has this property. However ${\mathcal R}$ is not
Cohen-Macaulay in that example; so it is natural to ask whether
the Cohen-Macaulayness of ${\mathcal R}$ implies the one of
${\mathcal F}$. Also this question has a negative answer:

\begin{Example}\label{example 3.8}
{\rm Let $R$ and $I$ be as in Example~\ref{example3.7}. By
adjoining two power series variables $x$ and $y$ we obtain the
ideal $I'=(I,x,y)\subset R[\![x,y]\!]$. Now $I'$ has height 3,
analytic spread 3, reduction number 2, and ${\mathcal G}(I')$ and
${\mathcal F}(I')$ are polynomial rings over ${\mathcal G}(I)$ and
${\mathcal F}(I)$, respectively. Thus ${\mathcal G}(I')$ is
Cohen-Macaulay and hence ${\mathcal R}(I')$ is Cohen-Macaulay by
\cite[3.6]{SUV}. However, ${\mathcal F}(I')$ still fails to be
Cohen-Macaulay.}
\end{Example}

Next, we list more cases in which the Cohen-Macaulayness of
${\mathcal R}$ and ${\mathcal G}$ implies the one of ${\mathcal
R}/K{\mathcal R}$.

\begin{Corollary}\label{corollary3.12}
Let $(R, {\mathfrak m})$ be a local Cohen-Macaulay ring with
infinite resi\-due field, and let $I$ be a perfect $R$-ideal of
height $2$ with analytic spread $\ell$ and minimal number of
generators $n$. Assume that $I$ satisfies $G_{\ell}$. If
${\mathcal R}$ is Cohen-Macaulay then ${\mathcal R}/K{\mathcal R}$
is Cohen-Macaulay for every ${\mathfrak m}$-primary ideal $K$
containing ${\rm Fitt}_{n-1}(I)$. In particular, if ${\mathcal R}$
is Cohen-Macaulay then ${\mathcal F}$ is Cohen-Macaulay.
\end{Corollary}
\demo We can apply Corollary~\ref{corollary3.9}, since $r \leq
\ell-1=\ell-g+1$ according to \cite[3.6]{SUV}. \QED

\medskip

\begin{Corollary}\label{cor3.6}
Let $(R, {\mathfrak m})$ be a local Gorenstein ring with infinite
residue field $k$, and let $I$ be an $R$-ideal of height $g$ with
analytic spread $\ell \geq g+1$ and minimal number of generators
$n=\ell+1$. Suppose that $I$ satisfies $G_{\ell}$ and that ${\rm
depth} \, R/I^j \geq {\rm dim} \, R/I -j +1$ whenever $1 \leq j
\leq \ell-g$. Furthermore assume that the natural map $I \otimes_R
k \longrightarrow {\rm Fitt}_{n-1}(I) \otimes_R k$ is not
injective. If ${\mathcal G}$ is Cohen-Macaulay then ${\mathcal
R}/K{\mathcal R}$ is Cohen-Macaulay for every ${\mathfrak
m}$-primary ideal $K$ containing ${\rm Fitt}_{n-1}(I)$. In
particular, if ${\mathcal G}$ is Cohen-Macaulay then ${\mathcal
F}$ is Cohen-Macaulay.
\end{Corollary}
\demo By \cite[2.3 and 2.4]{CP}, the ideal $I$ satisfies
$S_j(I)\cong I^j$ for $1 \leq j \leq \ell-g+1$ and is strongly
Cohen-Macaulay. In particular ${\rm depth}\, R/I^j \geq {\rm
dim}\, R/I -j +1$ whenever $1 \leq j \leq \ell-g+1$ \cite[the
proofs of 5.1 and 4.6]{HSV}. Now \cite[2.1]{PU} shows that $r =
\ell-g+1$, and the assertion follows from
Corollary~\ref{corollary3.9}. \QED

\medskip

\begin{Corollary}\label{cor3.7}
Let $(R, {\mathfrak m})$ be a local Gorenstein ring with infinite
residue field, and let $I$ be a perfect Gorenstein $R$-ideal of
height $3$ with analytic spread $\ell$ and minimal number of
generators $n$. Assume that $I$ satisfies $G_{\ell}$. If
${\mathcal G}$ is Cohen-Macaulay, then ${\mathcal R}/K{\mathcal
R}$ is Cohen-Macaulay for every ${\mathfrak m}$-primary ideal $K$
containing ${\rm Fitt}_{n-1}(I)$. In particular, if ${\mathcal G}$
is Cohen-Macaulay then ${\mathcal F}$ is Cohen-Macaulay.
\end{Corollary}
\demo We can apply Corollary~\ref{corollary3.9}, since $r \leq
\ell-2 =\ell-g+1$ according to \cite[3.1]{PU}. \QED

\medskip

It is natural to raise the question of whether the
Cohen-Macaulayness of ${\mathcal F}$ implies the one of ${\mathcal
R}$ and ${\mathcal G}$. The answer is negative even in the case of
perfect ideals of height two satisfying $G_{\ell}$. It is easy to
build counterexamples for ideals with second analytic deviation
one, because in this case, if $I$ is generated by homogeneous
polynomials of the same degree in a power series ring over a
field, $\mathcal F$ is a hypersurface ring and so it is always
Cohen-Macaulay. However, ${\mathcal R}$ is not Cohen-Macaulay if
the {\it row condition} is not satisfied.

We recall the following result from \cite[5.4]{U2}: Let $R$ be a
local Gorenstein ring with infinite residue field, let $I$ be a
perfect $R$-ideal of height 2, with analytic spread $\ell$ and
reduction number $r$. Assume that $I$ satisfies $G_{\ell}$. The
following are equivalent:
\begin{itemize}
\item[$({\it a})$]
$\mathcal R$ is Cohen-Macaulay.

\item[$({\it b})$]
$r<\ell$ (in which case $r=0$ or $r=\ell-1$).

\item[$({\it c})$]
${\rm Fitt}_{\ell}(I) = {\rm Fitt}_0(I/J)$ for some ideal $J
\subset I$ with $\mu(J)=\ell$.
\end{itemize}
Condition $({\it c})$ above is usually referred to as the `row
condition.'

\medskip

We end this section with a characterization of the Cohen-Macaulay
property of ${\mathcal G}$. If $I$ is an $R$-ideal of height $g$
with analytic spread $d = {\rm dim}\, R$ and the `expected
reduction number' $\leq d-g+1$, then Theorem~\ref{3.8} recovers an
earlier result of Johnson and Ulrich \cite[3.1]{JU} on the
Cohen-Macaulayness of ${\mathcal G}$. Without any restriction on
the reduction number, Theorem~\ref{3.8} characterizes the
Cohen-Macaulay property of ${\mathcal G}$ in terms of certain
intersection conditions, whose necessity was already known by
previous work of Polini and Ulrich \cite[1.2]{PU} $($see also
\cite[the proof of 5.2]{AbHu}$)$. Theorem~\ref{3.8} is also a
generalization of a well known criterion of Valabrega and Valla
\cite[2.7]{VV} to non ${\mathfrak m}$-primary ideals. The spirit
of this result and the methods of the proof are similar to our
earlier ones $($see for instance Theorem~\ref{2.10}$)$.

\begin{Theorem}\label{3.8}
Let $R$ be a local Cohen-Macaulay ring of dimension $d$ with
infinite residue field and let $I$ be an $R$-ideal of height $g$.
Assume that $I$ satisfies property $G_d$ and ${\rm depth}\, R/I^j
\geq {\rm dim}\, R/I -j+1$ for $1 \leq j \leq d-g$. Let $J$ be a
reduction of $I$ generated by $d$ elements with ${\rm ht}\, J
\colon I \geq d$ and write $r=r_J(I)$. Then the following are
equivalent:
\begin{itemize}
\item[$({\it i})$]
${\mathcal G}$ is Cohen-Macaulay;

\item[$({\it ii})$]
$JI^{j-1} \cap I^{j+1}=JI^j$ whenever $d-g+1 \leq j \leq r-1$.
\end{itemize}
In particular, if $r \leq d-g+1$ then ${\mathcal G}$ is
Cohen-Macaulay.
\end{Theorem}
\demo By \cite[1.2$({\it b})$]{PU} one has that $({\it i})$
implies $({\it ii})$. To prove the converse notice that ${\rm
depth}\, R/I^j \geq {\rm dim}\, R/I -j + 1$ whenever $1 \leq j
\leq d-g+1$. Furthermore if condition $({\it ii})$ holds then it
holds for every $j \geq d-g+1$. Let $a_1, \ldots, a_d$ be general
elements in $J$ and set ${\mathfrak a}_i=(a_1, \ldots, a_i)$ for
$0 \leq i \leq d$. We use the convention $I^j=R$ for $j \leq 0$.

We first claim that for $j \geq i-g$,
\begin{eqnarray}
{\mathfrak a}_i I^{j-1} \cap I^{j+1} & = & {\mathfrak a}_iI^j
\quad
\begin{array}{l}{\rm whenever \ \ } 0 \leq i \leq d-1 \\
{\rm or \ \ }  j \geq d-g+1,
\end{array}
\label{eq4} \\
({\mathfrak a}_iI^j \colon a_{i+1}) \cap I^{j+1} & = & {\mathfrak
a}_iI^j \quad\,\,\, {\rm whenever} \quad 0 \leq i \leq d-1.
\label{eq5}
\end{eqnarray}
We prove $(\ref{eq4})$ and $(\ref{eq5})$ simultaneously by
induction on $j$. The assertions are clear if $j<0$ since then $i
\leq g-1$. Now the claims follow for $j \leq d-g$ by applying
Lemma~\ref{lemma2.7} twice with $s=d$, $t=d-g$ and with $s=d$,
$t=d-g+1$, respectively. Next, suppose that $j \geq d-g+1$. We use
decreasing induction on $i$. If $i=d$ then $(\ref{eq4})$ follows
from assumption $({\it ii})$, whereas $(\ref{eq5})$ is vacuous.
Let $i \leq d-1$. As to $(\ref{eq4})$ we have that
\begin{eqnarray*}
{\mathfrak a}_iI^{j-1} \cap I^{j+1} & = & {\mathfrak a}_iI^{j-1}
\cap {\mathfrak a}_{i+1}I^{j-1} \cap I^{j+1} \\
& = & {\mathfrak a}_iI^{j-1} \cap {\mathfrak a}_{i+1}I^j
\hspace{1.8cm} {\rm by \ induction \ on \ }
i \ {\rm in \ } (\ref{eq4}) \\
& = & {\mathfrak a}_iI^{j-1} \cap ({\mathfrak a}_iI^j + a_{i+1}I^j) \\
& = & {\mathfrak a}_iI^j + ({\mathfrak a}_iI^{j-1} \cap a_{i+1}I^j) \\
& = & {\mathfrak a}_iI^j + a_{i+1}(({\mathfrak a}_iI^{j-1}
\colon a_{i+1}) \cap I^j) \\
& = & {\mathfrak a}_i I^j + a_{i+1}{\mathfrak a}_iI^{j-1}
\hspace{1.5cm} {\rm by \ induction \ on \ } j \ {\rm in \ } (\ref{eq5}) \\
& = & {\mathfrak a}_i I^j.
\end{eqnarray*}
As to $(\ref{eq5})$ we have that
\begin{eqnarray*}
({\mathfrak a}_iI^j \colon a_{i+1}) \cap I^{j+1} & \subset &
({\mathfrak a}_iI^{j-1} \colon a_{i+1}) \cap I^j \cap I^{j+1} \\
& = & {\mathfrak a}_iI^{j-1} \cap I^{j+1} \hspace{1.4cm}{\rm by \
induction \ on} \ j \
{\rm in} \ (\ref{eq5}) \\
& = & {\mathfrak a}_i I^j \hspace{2.85cm} {\rm by \ } (\ref{eq4}).
\end{eqnarray*}
This completes the proof of $(\ref{eq4})$ and $(\ref{eq5})$.

Let $a_1', \ldots, a_d'$ denote the images of $a_1, \ldots, a_d$
in $[{\mathcal G}]_1 =I/I^2$. We claim that
\begin{equation}\label{eq6}
[(a_1', \ldots, a_i') :_{{\mathcal G}} a_{i+1}']_j=[(a_1', \ldots,
a_i')]_j
\end{equation}
whenever $0 \leq i \leq d-1$ and $j \geq i-g+1$. $($See also
\cite[proof of 2.8$({\it b})$]{JU}.$)$ We may assume $j \geq 0$.
Let $u \in [(a_1', \ldots, a_i') :_{{\mathcal G}} a_{i+1}']_j$.
Then $u=x+I^{j+1}$ for some $x \in I^j$, and we have
\begin{eqnarray*}
a_{i+1}x  \in ({\mathfrak a}_iI^j+I^{j+2}) \cap {\mathfrak
a}_{i+1}I^j & = & {\mathfrak a}_iI^j+ (I^{j+2} \cap
{\mathfrak a}_{i+1}I^j) \\
& = & {\mathfrak a}_i I^j + {\mathfrak a}_{i+1}I^{j+1}
\hspace{1cm} {\rm by}
\ (\ref{eq4})\\
& = & {\mathfrak a}_iI^j+a_{i+1}I^{j+1}.
\end{eqnarray*}
So $a_{i+1}(x-y) \in {\mathfrak a}_iI^j$ for some $y \in I^{j+1}$.
Since $u=x+I^{j+1}=x-y+I^{j+1}$, replacing $x$ by $x-y$ we may
assume that $x \in ({\mathfrak a}_iI^j \colon a_{i+1}) \cap I^j$.
However $({\mathfrak a}_iI^j \colon a_{i+1}) \cap I^j = {\mathfrak
a}_i I^{j-1}$ by $(\ref{eq5})$. Thus $x \in {\mathfrak
a}_iI^{j-1}$ which implies $u \in (a_1', \ldots, a_i')$. This
completes the proof of $(\ref{eq6})$.

Now the Cohen-Macaulayness of ${\mathcal G}$ follows from
$(\ref{eq6})$ and \cite[3.7]{JU} applied to $S={\mathcal G}$ with
$\ell =d$ . Notice that the assumptions on the depth of
$S/{\mathfrak b}_i$ in \cite[3.7]{JU} are satisfied because ${\rm
depth} \, R/I^j \geq {\rm dim}\, R/I - j+1$ for $1 \leq j \leq
d-g+1$ $($see the proof of \cite[3.13]{JU}$)$. \QED

\medskip

\section{When is ${\mathfrak m}I$ integrally closed?}

\noindent Let $(R, {\mathfrak m})$ be a Noetherian local ring and
let $I$ be an $R$-ideal. The {\it integral closure} of $I$ is the
ideal $\overline{I}$ of all elements of $R$ that satisfy an
equation of the form $X^n+c_1X^{n-1}+\cdots+c_{n-1}X+c_n=0$, where
$c_j\in I^j$. Equivalently, $\overline{I}$ is the largest ideal
having $I$ as a reduction. The ideal $I$ is called {\it integrally
closed} in case $\overline{I} = I$. If $\overline{I^j} = I^j$ for
every $j \geq 1$, we say that $I$ is {\it normal}.

We are interested in finding conditions that imply the integral
closedness of the ideal ${\mathfrak m}I$. As explained earlier in
the introduction, the motivation for this question originates in
Wiles' work on semistable curves. Our results were inspired by the
paper of H\"ubl and Huneke \cite[1.3 and 1.5]{HuHu}, however in
many cases they hold in the more general setting of an arbitrary
${\mathfrak m}$-primary ideal $K$ rather than ${\mathfrak m}$
itself. First, we give a slight generalization of
\cite[1.5]{HuHu}.

\begin{Proposition} \label{theorem4.1}
Let $(R, {\mathfrak m})$ be a normal local ring of dimension $d$,
let $I$ be a normal $R$-ideal with analytic spread $d$, and let
$K$ be an ${\mathfrak m}$-primary ideal. Assume that ${\mathcal
R}/K{\mathcal R}$ is unmixed. Then $KI^j$ is integrally closed for
every $j \geq 0$. In particular $K$ is integrally closed.
\end{Proposition}
\demo We may assume $d>0$. Notice that $K{\mathcal R}$ is a
divisorial ideal of the normal domain ${\mathcal R}$. Such an
ideal is necessarily integrally closed. Since $K{\mathcal R} =
\displaystyle \bigoplus_{j \geq 0} KI^j$ we conclude that the
$R$-ideals $KI^j$ are integrally closed as well. \QED

\medskip

\begin{Corollary}\label{4.2}
Let $(R, {\mathfrak m})$ be a normal local Cohen-Macaulay ring of
dimension $d$ with infinite residue field, and let $I$ be a normal
strongly Cohen-Macaulay $R$-ideal of height $g$ with analytic
spread $d$, minimal number of generators $n$ and reduction number
$\leq d-g+1$. Assume that $I$ satisfies $G_d$. Then $KI^j$ is
integrally closed for every $j \geq 0$ and every ${\mathfrak
m}$-primary ideal $K$ containing ${\rm Fitt}_{n-1}(I)$. In
particular $K$ is integrally closed.
\end{Corollary}
\demo By Corollary~\ref{corollary3.9}, ${\mathcal R}/K{\mathcal
R}$ is Cohen-Macaulay. Now the assertion follows from
Proposition~\ref{theorem4.1}. \QED

\bigskip

The following is an extension of \cite[1.3]{HuHu}.

\begin{Proposition}\label{proposition4.3}
Let $(R, {\mathfrak m})$ be a Noetherian local ring with ${\rm
depth}\, R >0$, and let $I$ be an $R$-ideal. Let $s, t$ be
integers $\geq 1$. Assume that $I^s$ and $I^{s+t}$ are integrally
closed, ${\rm depth}\, R/I^s=0$ and ${\mathcal F}$ has no
homogeneous relations in degrees $\leq s+t$. Then $I^t \cap
\overline{{\mathfrak m}I^t} = {\mathfrak m}I^t$. In particular if
$I^t$ is integrally closed then ${\mathfrak m}I^t$ is integrally
closed as well.
\end{Proposition}
\demo See the proof of \cite[1.3]{HuHu}. \QED

\bigskip

We conclude this section by describing two occurrences of the
integral closedness of ${\mathfrak m}I$ in the case of perfect
ideals of height two and perfect Gorenstein ideals of height
three.

\begin{Corollary}\label{perf2}
Let $(R, {\mathfrak m})$ be a normal local Cohen-Macaulay ring of
dimension $d$, and let $I$ be a normal perfect $R$-ideal of height
$2$ with minimal number of generators $d+1$. Assume that $I$
satisfies $G_d$. Then ${\mathfrak m}I$ is integrally closed.
\end{Corollary}
\demo We may assume that the residue field of $R$ is infinite.
Notice that $\ell(I)=d$ by \cite[3.1]{H2} and \cite[1.11]{U}. Let
$r$ be the reduction number of $I$. If $r \leq d-1$ then the
assertion follows from Corollary~\ref{4.2}. If $r \geq d$ then
${\mathcal F}={\mathcal R}/{\mathfrak m} {\mathcal R}$ has no
homogeneous relations in degrees $\leq d$ since
$\mu(I)=\ell(I)+1$. As ${\rm pd}\, S_{d-1}(I)=d-1$ by
\cite[Theorem 1$({\it b})$]{Wey}, one has ${\rm depth} \,
S_{d-1}(I) =1$. This implies that $S_{d-1}(I) \cong I^{d-1}$ by
\cite[the proofs of 5.1 and 4.6]{HSV}, and then ${\rm depth} \,
R/I^{d-1} = 0$. Now Proposition~\ref{proposition4.3} applies with
$s=d-1$ and $t=1$. \QED

\medskip

\begin{Corollary}\label{gor3}
Let $(R, {\mathfrak m})$ be a local Gorenstein ring of dimension
$d$ with $d$ odd, and let $I$ be a perfect Gorenstein $R$-ideal of
height $3$ with minimal number of generators $n \geq d$. Assume
that $I$ satisfies $G_d$ and $I$, $I^{d-2}$, $I^{d-1}$ are
integrally closed. Then ${\mathfrak m}I$ is integrally closed.
\end{Corollary}
\demo By \cite[2.5, 2.15, 4.3$({\it g})$, 5.9, 6.17$({\it c})({\it
i})$ and $({\it d})({\it ii})$]{KU} ${\rm pd} \, S_{d-2}(I) = d-1$
and $S_j(I) \cong I^j$ for $1 \leq j \leq d-1$. In particular,
${\rm depth} \, R/I^{d-2} = 0$ and ${\mathcal F}={\mathcal
R}/{\mathfrak m}{\mathcal R}$ has no homogeneous relations in
degrees $\leq d-1$. Now applying Proposition~\ref{proposition4.3}
with $s=d-2$ and $t=1$ we conclude that ${\mathfrak m}I$ is
integrally closed. \QED

\end{document}